\documentclass[oneside]{article}
\usepackage{amssymb,latexsym,amsmath,amsthm,bm}
\usepackage{fancyhdr}
\usepackage{color}
\usepackage{mathrsfs}

\newtheorem{thm}{Theorem}[section]
\newtheorem*{thm*}{Theorem}
\newtheorem{prop}{Proposition}[section]
\newtheorem{lemma}{Lemma}[section]
\newtheorem{conjecture}{Conjecture}[section]

\newcommand{\ep}{\varepsilon}

\newcommand{\con}{\equiv}

\newcommand{\ndiv}{\nmid}
\newcommand{\modd}[1]{\; ( \text{mod} \; #1)}

\newcommand{\maps}{\rightarrow}
\newcommand{\intersect}{\cap}

\newcommand{\al}{\alpha}
\newcommand{\be}{\beta}

\newcommand{\del}{\delta}
\newcommand{\Del}{\Delta}

\newcommand{\Sig}{\Sigma}

\newcommand{\Acal}{\mathcal{A}}

\newcommand{\Bcal}{\mathcal{B}}

\newcommand{\Gcal}{\mathcal{G}}

\newcommand{\Mcal}{\mathcal{M}}

\newcommand{\xbf}{{\bf x}}

\newcommand{\Fbb}{\mathbb{F}}

\newcommand{\Z}{\mathbb{Z}}

\newcommand{\Pscr}{\mathscr{P}}

\newcommand{\beq}{\begin{equation}}
\newcommand{\eeq}{\end{equation}}
\numberwithin{equation}{section}

	\newcommand{\xtra}[1]{}

\begin{document}

\title{Burgess bounds for short mixed character sums} 
\author{D. R. Heath-Brown\footnote{Mathematical Institute, Radcliffe Observatory Quarter, Woodstock Road, Oxford
OX2~6GG, \emph{rhb@maths.ox.ac.uk}} \\
	  L. B. Pierce\footnote{Hausdorff Center for Mathematics, 62 Endenicher Allee, 53115 Bonn, \emph{pierce@math.uni-bonn.de}} }
\date{}
\maketitle  

\begin{abstract}
This paper proves nontrivial bounds for short mixed character sums by
introducing estimates for Vinogradov's mean value theorem into a 
version of the Burgess method. 

\end{abstract}

\section{Introduction}
Let $\chi(n)$ be a non-principal character of modulus $q$, and consider the character sum 
\beq\label{SNH}
 S(N,H) = \sum_{N < n \leq N+H} \chi(n).
 \eeq
The classical P\'{o}lya-Vinogradov inequality provides the bound 
\[|S(N,H)| \ll q^{1/2} \log q,\] which is nontrivial only if the length $H$ of the character sum is longer than $q^{1/2+\ep}$. 
In a classic series of papers, Burgess \cite{Bur57}, \cite{Bur62B},
\cite{Bur63A}, \cite{Bur86} introduced a method for bounding short
character sums that results in the following well-known bound: for
$\chi$ a primitive multiplicative character to a prime modulus $q$ one has 
\beq\label{Burgess}
 |S(N,H) |\ll H^{1-\frac{1}{r}}q^{\frac{r+1}{4r^2}}\log q,
 \eeq
for any integer $r\geq 1$; moreover this bound is uniform in $N$. 
This provides a nontrivial estimate for $S(N,H)$ as soon as $H >
q^{1/4+\ep}$; more precisely if $H=q^{1/4+\kappa}$, then the Burgess
bound is of size $Hq^{-\del}$ with $\del \approx \kappa^2$.  Indeed
Burgess proved a similar bound for arbitrary moduli $q$ when $r\le 3$,
and for general cube-free moduli for all $r$.

Burgess bounds have found valuable applications in a range of settings, and it would be highly desirable to develop variations of the Burgess method for mixed character sums of the form
\[ \sum_{N < n \leq N+H} e_q(f_1(n) \overline{f_2(n)}) \chi(f_3(n) \overline{f_4(n)}),\]
for appropriate polynomials $f_1,\ldots, f_4$ and $e_q(t) = e^{2\pi i t/q}$.
However, it has proved difficult to handle sums involving $\chi$ evaluated at anything other than a linear function of $n$.

This paper will be concerned with the short mixed character sum
\beq\label{Sf}
 S(f;N,H) =   \sum_{N < n \leq N+H} e( f(n)) \chi(n),
 \eeq
for prime moduli $q$, where $f$ is a real-valued polynomial and $e(t)
= e^{2\pi i 
   t}$. Recall that at its heart, the Burgess method involves breaking
 the range of the summand $n$ into residue classes modulo an auxiliary
 prime $p$. One then averages over a set of such primes $p$, and it is
 crucial that the argument of the characters may be made independent
 of $p$ (although the range of summation may still depend on $p$). 
More explicitly, fix a prime $p \ndiv q$ and split the set of $n \in (N,N+H]$ into residue classes modulo $p$ by writing $n=aq +pm$ with $0 \leq a <p$ and $m \in (N',N'+H']$ with $N' = (N-aq)/p,$ $H'=H/p$. Then, for example, the multiplicative character sum (\ref{SNH}) may be written as
\[ S(N,H) = \sum_{0 \leq a<p} \;\;\sum_{N' < m \leq N' + H'}
\chi(aq+pm) = \chi(p)  \sum_{0 \leq a<p} \;\;\sum_{N' < m \leq N' +
  H'} \chi(m) ,\] 
so that after averaging over a set $\Pscr$ of primes,
\beq\label{classic}
 |S(N,H)| \leq \frac{1}{|\Pscr|} \sum_{p \in \Pscr} \sum_{0 \leq a<p} 
\left| \sum_{N' < m \leq N' + H'} \chi(m) \right|.
 \eeq
The Burgess argument then proceeds by manipulating the intervals of summation in order to reach a complete character sum that may be bounded (in most cases) by the Weil bound. 
This reveals a fundamental barrier quickly reached by a naive application of the Burgess method to the mixed character sum (\ref{Sf}): it is not trivial to make the argument of the polynomial $f$ independent of sufficiently many primes $p$, and without this independence, averaging over auxiliary primes as in (\ref{classic}) cannot proceed successfully.

For the case of $f$ linear, Burgess \cite{Bur88} proved that for $f (n) = an/q$ with $0<a<q$ and $q$ prime,
\beq\label{S_linear}
 |S(f;N,H) | \ll H^{1-\frac{1}{r}}q^{\frac{1}{4(r-1)}} (\log q)^2,
 \eeq
for any $r \geq 2$ and $0<N,H < q$; this was later extended in \cite{Bur89} to the case $r=3$ and $q$ an arbitrary positive integer.
A similar result was also proved by Friedlander and Iwaniec \cite{FriIwa93}, as a consequence of more general bounds for weighted multiplicative character sums. 

 In a 1995 paper, Enflo \cite{Enf95} reported a nontrivial bound for $S(f;N,H)$ for $f$ a real-valued polynomial of any degree $d$ and $H=q^{1/2}$, with $q$ prime.
His proof introduced the idea of using Weyl differencing $d$ times
before applying the Burgess method, thus stripping off the exponential
factor $e(f(n))$ entirely. This insight removes the problem of dependence on the auxiliary primes, and allows the Burgess method to proceed.
A careful analysis of Enflo's method gives the following result:
\begin{thm}\label{thm_Enflo}
Let $f$ be a real-valued polynomial of degree $d$ and $\chi$ a
non-principal character to a prime modulus $q$. Then for any $r \geq 1$
and $H<q^{\frac{3}{4} + \frac{1}{4r}}$ we have
\[ \sum_{N<n \leq N+H} e(f(n)) \chi(n) \ll_{r,d,\ep} H^{1-\frac{1}{2^dr}} q^{\frac{r+1}{2^{d+2}r^2}+\ep},\]
uniformly in $N$.
\end{thm}

As this result is surpassed by new methods, we do not give a proof here. 
Note that this recovers the original Burgess bound (\ref{Burgess}) in the case $d=0$, and for any $d$ it proves a nontrivial bound as long as $H > q^{1/4+\ep}$. Note also that it is clear that an upper bound on $H$ is required as soon as $d\geq 1$. For example, if $f(n)=n/q$ and $H=mq$ for some $m \geq 1$ then $S(f;N,H)=mG_q(\chi)$, where $G_q(\chi)$ is the Gauss sum. Then $|S(f;N,H)|=Hq^{-1/2}$ precisely, so it is not possible to attain a generic upper bound of the form $H^\al q^\be$ with $\al<1$ for arbitrary $H$.

More recently, Chang \cite{Cha10} introduced another idea that allows one to remove the dependence of $e(f(n))$ on the auxiliary primes $p$. Roughly speaking, the idea is to approximate $S(f;N,H)$ by $S(\tilde{f};N,H)$, where $\tilde{f}$ has real coefficients that are sufficiently close to those of $f$ but are independent of $p$. Chang's result improves on that of Enflo, proving that as soon as $H > q^{1/4+\kappa}$,
\beq\label{sum_Chang_ref}
 \sum_{0<n \leq H} e(f(n)) \chi(n) \ll Hq^{-\del},
 \eeq
 where 
 \beq\label{Chang_del}
 \del = \frac{\kappa^2}{4((d+1)^2+2)(1+2\kappa)}.
 \eeq
 (In fact Chang's results in \cite{Cha10} apply more generally to mixed character sums over $\mathbb{F}_{q^n}$ for any $n \geq 1$.) 
Chang furthermore proved in \cite{Cha13} a result for square-free $q$ that is similar to (\ref{sum_Chang_ref}), but with an additional factor $\tau(q)^{4(\log d)d^{-2}}$.

A  refinement of Chang's argument improves the result to: 
\begin{thm}\label{thm_Chang}
Let $f$ be a real-valued polynomial of degree $d \geq 0$ and $\chi$ a
non-principal character to a prime modulus $q$. Set
\begin{equation}\label{DD}
D:=\frac{d(d+1)}{2}.
\end{equation}
Then if $r \geq 1$
and $H<q^{\frac{1}{2} + \frac{1}{4r}}$ we have 
\[ \sum_{N<n \leq N+H} e(f(n))\chi(n)  \ll_{r,d} H^{1-\frac{1}{r}}
q^{\frac{r+1 + D}{4r^2}}(\log q)^2, 	 \]
	 uniformly in $N$.
	 \end{thm}
We shall use the notation (\ref{DD}) throughout the paper.

We do not claim Theorem \ref{thm_Chang} as substantially new; the small improvement is a consequence of approximating the coefficients of monomials in $f$ more accurately for higher degree monomials; Chang approximates the coefficients with the same accuracy for every degree. Supposing that the result of Theorem \ref{thm_Enflo} achieves its minimum at a value $r_0$, we may compare it to the result of Theorem \ref{thm_Chang} for $r=2^d r_0$, and see that Theorem \ref{thm_Chang} is as strong for $d=1,2$ and stronger than Theorem \ref{thm_Enflo} for $d \geq 3$. Additionally, note that for $H<q^{\frac{1}{2} + \frac{1}{4r}}$, the bound of Theorem \ref{thm_Chang} is nontrivial only if $r \geq 1 + D$. 

If $H=q^{\frac{1}{4}+\kappa}$ for some small $\kappa>0$, then Theorem \ref{thm_Chang} yields a nontrivial bound $Hq^{-\del}$ where $\del$ behaves approximately like 
\beq\label{delta_Chang_refine}
\del = \frac{\kappa^2}{D +1}, 
\eeq
for sufficiently small $\kappa$ and sufficiently large $d$, and hence is approximately a factor of $8$ better than (\ref{Chang_del}).
(See Section \ref{sec_opt_Chang} for details.)

The novelty of this paper appears in the following strategy: by
choosing the coefficients of $\tilde{f}$ according to a certain grid,
we are able to introduce a nontrivial auxiliary averaging that leads
to a bound involving the number $J_{r,d}(X)$ occuring in Vinogradov's
mean value theorem.  This is the number of solutions to the system of Diophantine equations given by
\[ x_1^m + \cdots + x_r^m = x_{r+1}^m + \cdots + x_{2r}^{m}, \qquad 1
\leq m \leq d,\] 
 where $d$ is the degree of $f$ and $1 \leq x_1,\ldots, x_{2r} \leq X$.
The celebrated new results of Wooley (most recently \cite{Woo14bX}
\cite{Woo14aX}) on Vinogradov's mean value theorem provide
exceptionally sharp bounds for $J_{r,d}(X)$ and lead to a significant
improvement on Theorem \ref{thm_Chang}.  

Let us recall the main conjecture in the setting of Vinogradov's mean value theorem:
\begin{conjecture}\label{conj_VMT}
 For every $r \geq 1, d\geq 1$ and $\ep>0$, 
\beq\label{main_J}
 J_{r,d}(X) \ll_{r,d,\ep} X^\ep (X^r + X^{2r-D}).
 \eeq
 \end{conjecture}
 
 Conditional on this bound for $J_{r,d}(X)$ we prove our main result:
 
 \begin{thm}\label{thm_Vin_conj}
Let $f$ be a real-valued polynomial of degree $d \geq 1$ and $\chi$ a non-principal character to a prime modulus $q$. Assume
Conjecture \ref{conj_VMT} holds.  Then 
for integers $r > D$ and $H<q^{ \frac{1}{2} + \frac{1}{4(r- D)}}$ we have
\beq\label{conj3}
 \sum_{N<n \leq N+H} e(f(n))\chi(n) \ll_{r,\ep}
	H^{1-\frac{1}{r}} q^{\frac{r+1 -D}{4r(r -D)} + \ep},
		 \eeq
		uniformly in $N$, for any $\ep>0$.
\end{thm}	

The method of proof for Theorem \ref{thm_Vin_conj} also yields character sum bounds (conditional on Conjecture \ref{conj_VMT}) in the range $r \leq D$, but it turns out that these bounds are no better than trivial.
Note that the $d=0$ case of (\ref{conj3}) would recover the classical Burgess bound (\ref{Burgess}). For fixed $d$, in the limit as $r \maps \infty$, the bound (\ref{conj3}) is nontrivial for $H \geq q^{1/4+ \ep}$. 
A direct comparison shows that (\ref{conj3}) matches Theorem
\ref{thm_Chang} when $r=D+1$ (though the admissible range for $H$ is longer), and is sharper as soon as $r >D+1$.

If $H=q^{\frac{1}{4}+\kappa}$ for some small $\kappa>0$, then Theorem
\ref{thm_Vin_conj} would yield a nontrivial bound $Hq^{-\del}$ where $\del$
behaves approximately like  
\beq\label{delta_Vin}
\del = \left( \frac{2\kappa}{1 + \sqrt{1+4D\kappa}} \right)^2.
\eeq
(See Section \ref{sec_delta_Vin} for details.) 
For any fixed $d$, as $\kappa \maps 0$, this behaves like 
\[ \del = \kappa^2,\]
which we note is independent of $d$, and is in fact as strong as the
original Burgess bound for multiplicative character sums.

Note that for $d=1,2$, the bound of Conjecture \ref{conj_VMT} holds true trivially, for all $r \geq 1$. Thus the following are immediate corollaries of  Theorem \ref{thm_Vin_conj}:
\begin{thm}\label{thm_Vinogradov_1}
Let $f$ be a linear real-valued polynomial and $\chi$ a non-principal
character to a prime modulus $q$. Then for $r \geq 2$ and $H<q^{
  \frac{1}{2} + \frac{1}{4(r- 1)}}$ we have
\[ \sum_{N<n \leq N+H} e(f(n)) \chi(n) \ll_{r,\ep} H^{1-\frac{1}{r}} q^{\frac{1}{4(r-1)}+\ep}  ,\]
uniformly in $N$, for any $\ep>0$.
\end{thm}
Note that this generalizes the result (\ref{S_linear}) since $f$ may
now be any real-valued linear polynomial.

\begin{thm}\label{thm_Vinogradov_2}
Let $f$ be a quadratic real-valued polynomial and $\chi$ a
non-principal character to a prime modulus $q$.  Then for $r \geq 4$
and $H<q^{ \frac{1}{2} + \frac{1}{4(r-3)}}$ we have
\[ \sum_{N<n \leq N+H} e(f(n)) \chi(n) \ll H^{1-\frac{1}{r}} q^{\frac{r-2}{4r(r - 3)} + \ep}, \]
uniformly in $N$, for any $\ep>0$.
\end{thm}

Recent breakthroughs of Wooley have provided very strong results toward Conjecture \ref{conj_VMT}. At the time of writing, the conjecture is now known to hold for all $r$ if $d=3$ and for  $r \geq d(d-1)$ when $d \geq 4$ (see \cite{Woo14aX}), and for 100\% of the critical interval $1 \leq r \leq D$ (see \cite{Woo14bX}). 
In our application, the results of Wooley for large $r$ make the following cases of Theorem \ref{thm_Vin_conj} unconditional.
\begin{thm}\label{thm_Vinogradov_3}
Let $f$ be a real-valued polynomial of degree $3$ and $\chi$ a
non-principal character to a prime modulus $q$. Then for $r \geq 7$
and $H<q^{ \frac{1}{2} + \frac{1}{4(r-6)}}$ we have
\[ \sum_{N<n \leq N+H} e(f(n)) \chi(n) \ll_{r,\ep} H^{1-\frac{1}{r}} q^{\frac{r-5}{4r(r - 6)} + \ep}, \]
uniformly in $N$, for any $\ep>0$.
\end{thm}

For $d \geq 4$, we have:

\begin{thm}\label{thm_Vinogradov_3_large}
Let $f$ be a real-valued polynomial of degree $d \geq 4$ and $\chi$ a non-principal character to a prime modulus $q$. 
Then for  $r \geq d(d-1)$ and $H<q^{ \frac{1}{2} + \frac{1}{4(r- D)}}$
we have
\[ \sum_{N<n \leq N+H} e(f(n)) \chi(n) \ll_{r,\ep} H^{1-\frac{1}{r}} q^{\frac{r+1 - D}{4r(r - D)} + \ep},\]
uniformly in $N$, for any $\ep>0$.
\end{thm}

Finally, in the intermediate range $D < r < d(d-1)$, we apply the so-called approximate main conjecture of \cite{Woo14bX}, which states that for all $d \geq 4$,
\[ J_{r,d}(X) \ll X^{\Delta_{r,d}}(X^r + X^{2r- D}) \]
 where $\Del_{r,d} = O(d)$ (see Theorem 1.5 of \cite{Woo14bX}). 
 This results in the following:

\begin{thm}\label{thm_Vinogradov_med}
Let $f$ be a real-valued polynomial of degree $d \geq 4$ and $\chi$ a non-principal character to a prime modulus $q$. 
Then for 
$D < r < d(d-1)$
 and $H<q^{ \frac{1}{2} + \frac{1}{4(r- D+\Delta)}}$ we have
\[ \sum_{N<n \leq N+H} e(f(n)) \chi(n) \ll_{r,\ep}  H^{1-1/r} q^{\frac{r+1 - D+2\Del}{4r(r - D + \Del)} + \ep},\]
where 
\[ \Del = \Del_{r,d} = O(d) \]
 is as specified in \cite{Woo14bX}.
\end{thm}	 
\medskip

We have stated these results in terms of polynomials $f(n)$. However
it is clear in principle that one can prove estimates for suitable general real-valued functions $f(n)$ by approximating them by appropriate
polynomials. Moreover, these methods can be extended to certain multi-variable sums.  We intend to return to this issue in the near future.

Although in this paper we shall confine ourselves to prime moduli $q$,
most of our results can be modified to apply to general square-free
moduli.  
In some cases however we cannot handle the full range 
$r>D$ occuring in Theorem \ref{thm_Vin_conj}.
We leave the details to the reader. 

For our proofs it will be convenient to assume that $d\ge 1$.  This
enables us to replace the use of the
Menchov-Rademacher device (originating in \cite{Men23}, \cite{Rad22})
by the simpler ``partial summation by Fourier series'' of Bombieri and
Iwaniec  \cite{BoIw86}.  Of course Theorem \ref{thm_Chang} remains true for
$d=0$, since it reduces to Burgess's bound (\ref{Burgess}).

\section{The Burgess method with coefficient
  approximation}\label{sec_Chang} 
To begin the proof of Theorems \ref{thm_Chang} and \ref{thm_Vin_conj},
we consider 
\[ T_d(N,H,\chi)=T(N,H) = \sup_{{\rm deg}(f)=d}\;\sup_{K\le H}
\left|\sum_{N < n \leq N+K} e(f(n)) \chi(n)\right|,\]
where $f$ runs over real-valued polynomials and $\chi$ is a non-principal
multiplicative character to a prime modulus $q$. We first note that
$T(N,H)$ has period $q$ with respect to $N$, so that we can assume
from now on that $0\le N<q$.

Fix a set of primes $\Pscr = \{ P < p \leq 2P \}$ for some parameter
$P\leq H$ that we will choose later.  Since $H=o(q)$ in all our
theorems we will have $p \ndiv q$ for $p\in\Pscr$. Hence we can split
$n \in (N,N+K]$ into 
  residue classes modulo $p$ by writing $n=aq +pm$ with $0 \leq a
  <p$.  This produces values
$m \in (N_{a,p},N_{a,p}+K_{a,p}]$ with $N_{a,p}= (N-aq)/p$ and
$K_{a,p}=K/p\le H/P$. 
Then 
\[\sum_{N < n \leq N+K} e(f(n)) \chi(n)
= \sum_{0 \leq a<p}\;\; \sum_{N_{a,p} < m \leq N_{a,p}+ K_{a,p}} e(f(aq+pm))\chi(aq+pm),\]
and as a result
\[ T(N,H) \leq \sum_{0 \leq a<p} T(N_{a,p},H/P).\]
We proceed to average over $\Pscr$, producing
\beq\label{Sap_sum}
T(N,H) \leq |\Pscr|^{-1}\sum_{p\in\Pscr}\;\sum_{0 \leq a<p} T(N_{a,p},H/P).
\eeq

We now use the following lemma.
\begin{lemma}
For any real number $L\ge 1$ we have
\beq\label{TOM}
T(U,L) \leq 4L^{-1} \sum_{U-L< m \leq U}T(m,2L).
 \eeq
 \end{lemma}
 
To see this, note that 
\[T(U,L)=\left|\sum_{U < n \leq U+K} e(f(n)) \chi(n)\right|\]
for some polynomial $f$ and some positive real number $K\le L$.  
Moreover
if $U-L<m\le U$ then
\[\sum_{U< n\leq U+K}e(f(n))\chi(n)=\sum_{m<n\leq U+K}e(f(n))\chi(n)
-\sum_{m< n\leq U}e(f(n))\chi(n),\]
whence
\[\left|\sum_{U< n\leq U+K}e(f(n))\chi(n)\right|\le 2T(m,2L),\]
since $U+K\le m+2L$.  The result then follows since the interval
$(U-L,U]$ contains at least $L/2$ integers $m$.

Applying (\ref{TOM}) to (\ref{Sap_sum}) with $U=N_{a,p}$ and $L=H/P$, 
we may conclude that
\begin{eqnarray*}
T(N,H)&\ll& |\Pscr|^{-1} (H/P)^{-1} \sum_{p\in\Pscr}\;\;\sum_{0 \leq a<p}
\;\;\sum_{N_{a,p}-H/P<m\le N_{a,p}}T(m,2H/P)\\
&\ll& H^{-1}(\log q) \sum_{p\in\Pscr}\;\;\sum_{0 \leq a<p}
\;\;\sum_{N_{a,p}-H/P<m\le N_{a,p}}T(m,2H/P),
\end{eqnarray*}
on noting that $|\Pscr|\gg P(\log P)^{-1}\gg P(\log q)^{-1}$.
We now define 
\[\Acal(m)=\#\left\{(a,p):\frac{N-aq}{p}-\frac{H}{P}<m\leq\frac{N-aq}{p}
\right\},\]
which allows us to write
\beq\label{SSS}
T(N,H)\ll 
H^{-1}(\log q)\sum_{m\in\mathbb{Z}} \Acal(m)T(m,2H/P).
\eeq

We then set
\[S_1 =  \sum_m \Acal(m)\]
and
\[S_2 = \sum_m \Acal(m)^2,\]
and we note the following facts, which we will prove in Section
\ref{sec_lemma_A'}.
\begin{lemma}\label{lemma_A'}
We have $\Acal(m)=0$ unless $|m| \leq 2q$. 
Moreover if $HP<q$ then $S_1\le S_2 \ll HP$.
\end{lemma}

From a repeated application of H\"{o}lder's inequality, it then
follows from (\ref{SSS}) that
\begin{eqnarray*}
T(N,H)&\ll& H^{-1}(\log q)S_1^{1-\frac{1}{r}}S_2^{\frac{1}{2r}}
\left\{\sum_{|m| \leq 2q}T(m,2H/P)^{2r}\right\}^{\frac{1}{2r}}\\
&\ll&
H^{-\frac{1}{2r}}P^{1-\frac{1}{2r}}(\log q)
\left\{\sum_{|m| \leq 2q}T(m,2H/P)^{2r}\right\}^{\frac{1}{2r}}.
\end{eqnarray*}
As previously noted, the function $T(m,K)$ is periodic in $m$, with
period $q$, so that in fact we have
\beq\label{TB}
T(N,H)\ll H^{-\frac{1}{2r}}P^{1-\frac{1}{2r}}(\log q)
\left\{\sum_{m=1}^qT(m,2H/P)^{2r}\right\}^{\frac{1}{2r}}.
\eeq

For any $M$ and $K>0$ we now define 
\[T_0(M,K)=\sup_{{\rm deg}(f)=d}\left|\sum_{M< n \leq M+K}
  e(f(n))\chi(n)\right|.\]
We can relate $T(M,K)$ to $T_0(M,K)$ using the following lemma, which
is an immediate consequence of Lemma 2.2 of Bombieri and Iwaniec \cite{BoIw86}.
\begin{lemma}
Let $a_n$ be a sequence of complex numbers supported on the integers
$n\in(A,A+B]$, and let $I$ be any subinterval of $(A,A+B]$.  Then
\[\sum_{n\in I}a_n\ll \big(\log(B+2)\big)\sup_{\theta\in\mathbb{R}}
\left|\sum_{A<n\le A+B}a_ne(\theta n)\right|.\]
\end{lemma}
Thus if  $d\ge 1$ and $K\le q$ then
\[T(M,K)\ll  T_0(M,K)\log(K+2)\ll T_0(M,K)\log q.\]
This is the only place in the argument where the condition $d\ge 1$ is used.
We now see that (\ref{TB}) becomes
\beq\label{SPHT}
T(N,H)\ll H^{-\frac{1}{2r}}P^{1-\frac{1}{2r}}(\log q)^2
S_3(2H/P)^{\frac{1}{2r}},
\eeq
where we have set
\[S_3(K)=\sum_{m=1}^qT_0(m,K)^{2r}.\]

We proceed to develop a bound for $S_3(K)$, under the assumption that
$K\le q$.
Having removed the maximum over the length of our intervals we now
handle the maximum over the polynomials $f$.  In effect we do this by
replacing the maximum by a sum over all ``distinct'' polynomials
modulo 1. The principle here is that two polynomials will be
effectively equivalent if their coefficients are sufficiently close.

Let $Q\geq K$ be an integer parameter to be chosen in due course.
We partition $[0,1]^{d+1}$ into boxes $B_\al$ of side-length $Q^{-j}$ in the
$j$-th coordinate, for $j=0,\ldots, d$. Note that the total number of
boxes is $Q^{D}$. For each box $B_\al$, fix $\theta_\al = (\theta_{\al,0},\ldots,
\theta_{\al,d})$ to be the vertex of $B_\al$ with the least value in
each coordinate. Thus each $\theta_{\al}$ takes the form  
\[(c_0Q^{-0}, c_1 Q^{-1}, \ldots, c_d Q^{-d})\]
for some integers $0 \leq c_j \leq Q^j-1$, $0 \leq j \leq d$. 
(Chang's original argument \cite{Cha10} chooses the boxes to be of side-length
$Q^{-d}$ in all coordinates, and allows $\theta_\al$ to be any point
in the box $B_\al$.) 
Define for any $\theta \in [0,1]^{d+1}$ the polynomial 
\[ \theta(X) := \sum_{j=0}^d \theta_j X^j.\]
For any integer $m$, positive real number $t$, and index $\alpha$, set
\[T(\alpha;m,t):=\left|\sum_{0<n\le t} e(\theta_{\alpha}(n)) \chi(n+m)\right|.\]
We use these sums to approximate $T_0(m,K)$ as follows.

\begin{lemma}\label{lemma_tab}
Given an integer $m$ and real numbers $Q\ge K>0$, there is an 
index $\alpha$ such that
\[T_0(m,K)\ll_d T(\alpha;m,K)+K^{-1}\int_0^K T(\alpha;m,t)dt.\]
\end{lemma}
To prove this we observe that for integral $m$ we have
\begin{eqnarray*}
T_0(m,K)&=&\sup_{{\rm deg}(f)=d}\left|\sum_{m< n \leq m+K}
  e(f(n))\chi(n)\right|\\
&=&\sup_{{\rm deg}(f)=d}\left|\sum_{0< n \leq K}
  e(f(n))\chi(n+m)\right|.
\end{eqnarray*}
Suppose then that
\[T_0(m,K)=\left|\sum_{0<n \leq K} e(f(n))\chi(n+m)\right|\]
for some polynomial $f$ of degree $d$,
and write $f(X)=f_dX^d+\ldots+f_0$.  Clearly we may assume that
$0\le f_j\le 1$ for $0\le j\le d$. We then choose $\alpha$ so that
$|f_j-\theta_{\alpha,j}|\le Q^{-j}$ for each index $j$ and temporarily
write $\delta_j=f_j-\theta_{\alpha,j}$ for notational convenience.
Then, by summation by parts, we have
\begin{eqnarray*}
&& \hspace{-2cm} \sum_{0<n \leq K} e(f(n)) \chi(n+m)\\
\hspace{1cm} 
&=&\sum_{n \leq K} e\left(\sum_{j=0}^d \delta_j n^j\right)
  e(\theta_{\alpha}(n)) \chi(n+m)\\
&=& e\left(\sum_{j=0}^d \delta_j K^j\right)
\sum_{n \leq K} e(\theta_{\alpha}(n)) \chi(n+m)\\
&&\hspace{1cm}\mbox{} 
-\int_0^K \left\{\sum_{n \leq t} e(\theta_{\alpha}(n)) \chi(n+m)\right\}
\frac{d}{dt}  e\left(\sum_{j=0}^d\delta_j t^j\right) dt.\\ 
	\end{eqnarray*}
Since $|\delta_j|\le Q^{-j}$ we have 
\[ \left| \frac{d}{dt}  e\left(\sum_{j=0}^d \delta_j t^j\right)
\right| \leq 2\pi \sum_{j=1}^d j |\delta_j| t^{j-1}  \leq 2\pi 
\sum_{j=1}^d jQ^{-j} K^{j-1},\]
for $0\le t \leq K$.
Thus if $Q \gg K$ we have
\[ \left| \frac{d}{dt}  e\left(\sum_{j=0}^d \delta_j t^j\right)
\right| \ll_d  K^{-1}\]
and hence
\[\sum_{n \leq K} e(f(n))\chi(n+m)\ll_d T(\alpha;N,K) + 
K^{-1} \int_0^K T(\alpha;N,t)dt,\]
which proves the lemma.

An application of H\"older's now allows us to deduce from Lemma
\ref{lemma_tab} that
\[T_0(m,K)^{2r}\ll_d T(\alpha;m,K)^{2r}+
K^{-1}\int_0^K T(\alpha;m,t)^{2r}dt\]
for some index $\alpha$ depending on $m$ and $K$. This
dependence is rather awkward, and we circumvent it in the most trivial
way by summing over all available indices $\alpha$, giving
\[T_0(m,K)^{2r}\ll_d \sum_{\alpha}T(\alpha;m,K)^{2r}+
K^{-1}\sum_{\alpha}\int_0^K T(\alpha;m,t)^{2r}dt.\]
Thus
\beq\label{S2d5'}
S_3(K)\ll_d S_4(K)+K^{-1}\int_0^K S_4(t)dt
\eeq
if $0<K\ll Q$, where we have defined
\[S_4(\tau)=\sum_{\alpha}\sum_{m=1}^q T(\alpha;m,\tau)^{2r}.\]

Thus we now turn our attention to bounding the sum $S_4(\tau)$.
Recall the definition of the boxes $B_\al$,
and in particular the definition of the vertices $\theta_\al$. 
If $\xbf=(x_1,\ldots,x_{2r})$ we write
\[\Sigma_A(\xbf;q)=\sum_\al e\left(\sum_{i=1}^{2r} \ep(i) \theta_\al(x_i)\right),\]
where $\ep(i) =(-1)^i$. We also set
\[\Sigma_B(\xbf;\chi,q)=\sum_{m=1}^q \chi(F_\xbf(m))\]
where the polynomial $F_\xbf(X)$ is defined by
\beq\label{FxX}
 F_\xbf(X)  = \prod_{i=1}^{2r} (X+ x_i)^{\del_q(i)}.
 \eeq
Here $\del_q(i)=1$ if $i$ is even and $=\Del(q)-1$ if $i$ is odd,
where $\Del(q)$ is the order of the character $\chi$ modulo $q$.  

With this notation we then see upon expanding the sum that
\beq\label{al_m}
S_4(\tau)=  \sum_\al \sum_{m=1}^q T (\al; m,\tau)^{2r}
	= \sum_{\substack{\xbf\\ 0 < x_i \leq \tau}} 
\Sigma_A(\xbf;q) \Sigma_B(\xbf;\chi,q).
	\eeq

We will first prove Theorem \ref{thm_Chang} by averaging trivially over the boxes $B_\al$ and running the Weil bound argument that is typically found in applications of the Burgess method. The key proposition for Theorem \ref{thm_Chang} is:
\begin{prop}\label{prop_S4_Chang}
Suppose $q$ is prime.  Then for any $\tau\le q$ we have
\beq\label{T2r_Chang}
S_4(\tau)= \sum_\al \sum_{m=1}^q T (\al; m,\tau)^{2r} \ll_r 
Q^{D} (\tau^r q + \tau^{2r} q^{1/2}).
 \eeq
\end{prop}

Second, we will improve on this by averaging nontrivially over the
boxes $B_\al$, resulting in the key proposition for Theorem
\ref{thm_Vin_conj}: 
\begin{prop}\label{prop_S4}
Suppose $q$ is prime. Then for any $\tau\le q$ we have
\beq\label{T2r}
S_4(\tau)= \sum_\al \sum_{m=1}^q T (\al; m,\tau)^{2r} \ll_r 
Q^D  \left(\tau^r q+  J_{r,d}(\tau)q^{1/2} \right).  
 \eeq
\end{prop}

The propositions will be proved and the resulting theorems deduced in  Sections \ref{sec_mult_comp}
and \ref{s4}, respectively. Although Proposition
\ref{prop_S4_Chang} is an immediate consequence of Proposition
\ref{prop_S4} we have chosen to state and prove Proposition
\ref{prop_S4_Chang} separately, in order to highlight the different
aspects of our treatment.

\section{The multiplicative component}\label{sec_mult_comp}
We first consider the multiplicative character sum $\Sigma_B(\xbf;\chi,q)$. 
The well-known Weil bound implies the following:
\begin{lemma}\label{lemma_Weil}
Let $\chi$ be a character of order $\Del(q) >1$ modulo a prime
$q$. Suppose that $F (X)$ is a polynomial which is not a perfect
$\Del(q)$-th power over $\overline{\mathbb{F}}_q[X]$. Then
\[  \left| \sum_{m=1}^q \chi(F(m))\right| \leq ({\rm deg}(F)-1)\sqrt{q}. \]
\end{lemma}
We can apply Lemma \ref{lemma_Weil} to show that
$\Sigma_{B}(\xbf;\chi,q)$ is bounded by $O_{r}(q^{1/2})$, unless
the polynomial $ F_\xbf (X)$ is a perfect $\Del(q)$-th power over
$\overline{\Fbb}_{q}$. 
We define $\xbf = (x_1, \ldots, x_{2r})$ to be bad if for all
$i=1\ldots, 2r$, there exists $j \neq i$ such that $x_j=x_i$, and
$\xbf$ to be good otherwise. We take $\Bcal(\tau)$ to be the collection of
bad $\xbf$ with $0< x_i \leq \tau$ and similarly $\Gcal(\tau)$ to be the
collection of good $\xbf$ with $0<x_i \leq \tau$. The following is
immediate: 
\begin{lemma}\label{lemma_bad}
There are at most $r^{2r+1}\tau^r$ bad $\xbf$ with $0<x_i \leq \tau$, so that 
\beq\label{Bad_num}
\# \Bcal(\tau) \ll_r \tau^r.
\eeq
\end{lemma}
For the proof of the lemma we write the set $\{x_1,\ldots,x_{2r}\}$ without
repetitions as $\{y_1,\ldots,y_t\}$, say, where $t\le r$ since $\xbf$ is
bad.  We may suppose that the $y_i$ are arranged in ascending order.
There are at most $rK^r$ 
choices for such a set $\{y_1,\ldots,y_t\}$, and at most $r^{2r}$ choices
for $\xbf$ which correspond to each such set.  This suffices for the lemma.

Furthermore:
\begin{lemma}\label{lemma_A1k}
Fix $\xbf$ with $0<x_i \leq \tau$ for each $i=1,\ldots, 2r$ and fix a
prime $q$.  If $\tau\le q$ and $F_\xbf(X)$ is a perfect $\Del(q)$-th
power modulo $q$, then $\xbf$ is bad. 
\end{lemma}
This is obvious since if there were only one index $i$ for which $x_i$
takes a given value $y$ say, then the factor $X+y$ occurs in
$F_\xbf(X)$ with multiplicity either 1 or $\Del(q)-1$, neither of
which is divisible by $\Del(q)$.

If $\xbf$ is bad, we will apply the trivial bound $O(q)$ to $\Sigma_B(\xbf;\chi,q)$; we may conclude from (\ref{Bad_num}) that 
\beq\label{sum_bad}
 \sum_{\xbf \in \Bcal(\tau)} \left| \sum_{m=1}^q \chi(F_\xbf(m))\right| 
\ll_r \tau^r q.
 \eeq

For good $\xbf$ we may apply Lemmas \ref{lemma_Weil} and
\ref{lemma_A1k} to obtain the following standard result.
\begin{lemma}\label{lemma_sf}
If $q$ is prime and $\tau\le q$ then
\beq\label{sum_good}
\sum_{\xbf \in \Gcal(\tau)} | \sum_{m=1}^q \chi(F_\xbf(m))|  \ll_r 
\tau^{2r}q^{1/2}.
\eeq
\end{lemma}

\subsection{Proof of Theorem \ref{thm_Chang}}\label{sec_thm_Chang}

At this point we may prove Proposition \ref{prop_S4_Chang}. Using the trivial bound 
\[| \Sig_A(\xbf;q)| \leq Q^D\]
in (\ref{al_m}), we observe that
\begin{eqnarray*}
 \lefteqn{ \sum_\al \sum_{m=1}^q T(\al;m,\tau)^{2r}}\\ 
  	&\leq & Q^D  \left( \sum_{\xbf \in \Gcal(\tau)} \left| 
\sum_{m=1}^q \chi(F_\xbf(m) )\right| +\sum_{\xbf \in \Bcal(\tau)} \left| 
\sum_{m=1}^q \chi( F_\xbf(m) ) \right|\right) .
	\end{eqnarray*}
We substitute the bounds (\ref{sum_good}) and (\ref{sum_bad}) to
complete the proof of Proposition \ref{prop_S4_Chang}. Applying
Proposition \ref{prop_S4_Chang} to $S_4(K)$ and $S_4(t)$ in
(\ref{S2d5'}), we may
conclude that for any $K \le q$ we have
\[S_3(K) \ll_{r,d} Q^D (K^{2r}q^{1/2} + K^r q)\]
so long as the integer $Q$ is at least $K$. We apply this in
(\ref{SPHT}) with $K=2H/P$ and $Q=\lceil 2H/P\rceil$, obtaining
\[ T(N,H) \ll_{r,d} H^{-\frac{1}{2r}} P^{1-\frac{1}{2r}}(\log q)^2 
(H/P)^{\frac{D}{2r}}  \left((H/P)^{2r}q^{1/2} + (H/P)^r
  q\right)^{\frac{1}{2r}}. \]
We then extract the best result by choosing $P$ such that 
\[\frac{1}{2}Hq^{-1/(2r)} \leq P\leq Hq^{-1/(2r)}. \]
The  restriction $HP < q$ of Lemma
\ref{lemma_A'} is then satisfied when $H<q^{\frac{1}{2} + \frac{1}{4r}}$, and
we will also have $2H/P\le q$ for sufficiently large $q$. We therefore
obtain the result of Theorem \ref{thm_Chang} in the form
\[T(N,H) \ll_{r,d} H^{1-\frac{1}{r}}q^{\frac{r+1 + D}{4r^2}}(\log q)^2.\]

\subsection{Optimal choice of $r$}\label{sec_opt_Chang}
Recall that we have set
\[ D = \frac{1}{2} d(d+1).\]
We observe that if $H=q^{\frac{1}{4} + \kappa}$ for small $\kappa >0$, then the bound of Theorem \ref{thm_Chang} is of the form $Hq^{-\del}$ where 
\[ \del = \frac{\kappa r - \frac{1}{4} (D +1)}{r^2}.\]
As a function of $r$, this attains a maximum at the real value
\[ r(\kappa,d) := \frac{\frac{1}{2}(D+1)}{\kappa}.\]
Upon choosing the closest integer $r = r(\kappa,d) + \theta$ where
$-1/2 < \theta \leq 1/2$, we compute that for this choice of $r$ we have
\[ \del = \kappa^2 \left( \frac{ \frac{1}{4}(D+1) + \kappa
    \theta}{\frac{1}{4}(D+1)^2 + (d+1)\kappa \theta + \kappa^2
    \theta^2} \right).\] 
 For sufficiently small $\kappa$ this behaves like
\[ \del = \frac{\kappa^2}{D+1} .\]

\xtra{Note: finding optimal $r$ is equivalent to finding the maximum in $x$ of the function $f(x) = \frac{ax-b}{x^2}$ where $a=\kappa$, $b = \frac{1}{4}(\frac{1}{2}d(d+1) +1)$. Then $f'(x)=0$ when $x=2b/a$, which leads to our choice of $r$.}

\section{Introduction of the Vinogradov bounds}\label{s4}
We improve on the strategy of Theorem \ref{thm_Chang} by treating the
additive character sum $\Sigma_A(\xbf;q)$ in (\ref{al_m}) nontrivially. 
Recalling the definition of the vector $\theta_\al = (\theta_{\al,1}, \theta_{\al,2}, \ldots, \theta_{\al,d})$, we see that 
\begin{eqnarray*}
 \sum_\al e\left(\sum_{i=1}^{2r} \ep(i) \theta_\al(x_i)\right) 
	&= &\sum_\al e \left( \theta_{\al,1}\sum_{i=1}^{2r} \ep(i) x_i + \cdots +  \theta_{\al,d}\sum_{i=1}^{2r} \ep(i) x_i^d\right)\\
	& = & \prod_{s=1}^d  \left( \sum_{c=1}^{Q^s} e \left( \frac{c\sum_{i=1}^{2r} \ep(i) x_i^s}{Q^s} \right) \right)\\
		& = & Q^D\Xi_Q(\xbf),
	\end{eqnarray*}
say, where $\Xi_Q(\xbf)$ is the indicator function for the set 
\[  \{ \xbf = (x_1,\ldots, x_{2r}) \in \mathbb{N}^{2r}\cap(0,\tau]^{2r}:
\sum_{i=1}^{2r} \ep(i) x_i^s \con 0 \modd{Q^s}, \; \forall s\le d\}.\]

Our application has $0\le\tau\le K$ in (\ref{S2d5'}), and $Q\ge
K$ in Lemma \ref{lemma_tab}.  Moreover we will be taking $K=2H/P$ in
(\ref{SPHT}). Any integer $Q\ge 2H/P$ is therefore acceptable.  In the
definition of $\Xi_Q(\xbf)$ we will have
\[\left|\sum_{i=1}^{2r} \ep(i) x_i^s \right|<2r\tau^s\le 2rK^s\le
(2rK)^s=(4rH/P)^s.\]
Thus, by taking $Q = \lceil
4rH/P\rceil$, the congruences in the set above can hold only if they
are actually equalities in $\Z$.  We may then replace $\Xi_Q(\xbf)$ by the indicator function $\Xi(\xbf)$ of the set
\[  V_{r,d}(\tau) := \{ \xbf = (x_1,\ldots, x_{2r}) \in
\mathbb{N}^{2r}\cap(0,\tau]^{2r}: \sum_{i=1}^{2r} \ep(i) x_i^s =0, \; \forall
s\le d \}.\]

Then we see that (\ref{al_m}) may be bounded by
\[  \sum_\al \sum_{m=1}^q T(\al;m,\tau)^{2r} \leq
Q^D\{\Sigma(\Gcal)+\Sigma(\Bcal)\},\]
where
\[\Sigma(\Gcal)=\sum_{\xbf \in \Gcal (\tau)\intersect V_{r,d}(\tau)} 
\left| \sum_{m=1}^q \chi(F_\xbf(m) )\right|\]
and
\[\Sigma(\Bcal)=\sum_{\xbf \in \Bcal (\tau)\intersect V_{r,d}(\tau)} 
\left| \sum_{m=1}^q \chi(F_\xbf(m) )\right|.\]
We now prove Proposition \ref{prop_S4}.
Lemma \ref{lemma_A1k} shows that
$F_{\xbf}(X)$ is not a perfect $\Delta(q)$-th power modulo $q$ for 
$\xbf \in \Gcal(\tau)$ and $\tau\le q$, and
then Lemma \ref{lemma_Weil} yields
\[\sum_{m=1}^q \chi(F_\xbf(m) ) \ll_r q^{1/2}.\]
We expect $\xbf$ to be good generically, so we will apply the upper bound 
\[\#(\Gcal(\tau) \intersect V_{r,d}(\tau) )\leq \# V_{r,d}(\tau) = J_{r,d}(\tau),\]
whence
\beq\label{GV_prime}
\Sigma(\Gcal)=\sum_{\xbf \in \Gcal(\tau) \intersect V_{r,d}(\tau)} 
\left| \sum_{m=1}^q \chi(F_\xbf(m) )\right| \ll_r J_{r,d}(\tau)q^{1/2}.
\eeq
For $\xbf \in \Bcal(K)$ we use (\ref{sum_bad}) to deduce that
\[\Sigma(\Bcal)=\sum_{\xbf \in \Bcal (\tau)\intersect V_{r,d}(\tau)} 
\left| \sum_{m=1}^q \chi(F_\xbf(m) )\right|\leq
\sum_{\xbf \in \Bcal (\tau)} 
\left| \sum_{m=1}^q \chi(F_\xbf(m) )\right|\ll_r \tau^r q.\]
Proposition \ref{prop_S4} then follows.

\subsection{Proof of Theorem \ref{thm_Vin_conj}}\label{sec_thm_Vin}
We proceed to prove Theorem \ref{thm_Vin_conj}. Assuming that
Conjecture \ref{conj_VMT} holds, we see from Proposition \ref{prop_S4}
that  
\beq\label{T_conj}
S_4(\tau)\ll_{r,d,\ep} Q^D  \left\{(\tau^{r}
   + \tau^{2r-D})q^{1/2} + \tau^r q\right\}q^\ep . 
 \eeq
If $r \leq D$, the contribution of bad $\xbf$ dominates, and we cannot
obtain a nontrivial bound. Thus from now on we only consider $r> D$.
Since $d$ is then bounded in terms of $r$, the implied constant in the
$\ll_{r,d,\ep}$ notation may be bounded as a function of $r$ and $\ep$
alone.  We now
apply (\ref{T_conj}) to (\ref{S2d5'}) to conclude that for any $1
\leq K \le q$ we have
\[S_3(K) \ll_{r,\ep} Q^D(K^{2r-D}q^{1/2} + K^r q)q^{\ep}.\]
 We  apply this to (\ref{SPHT}) to obtain
\[ T(N,H) \ll_{r,\ep} H^{-\frac{1}{2r}}P^{1-\frac{1}{2r}}Q^{\frac{D}{2r}}
\left(K^{2r-D}q^{1/2}+K^rq\right)^{\frac{1}{2r}}q^{\ep}.\]
As before we take $K=2H/P$ and $Q=\lceil 4rH/P\rceil$.
  It is optimal to choose $P$ to balance the last two terms by taking
  \beq\label{P_choice}
\frac{1}{2} Hq^{- \frac{1}{2(r-D)}} \leq P<Hq^{- \frac{1}{2(r-D)}}.
  \eeq
We may then satisfy the requirement $HP<q$ of Lemma \ref{lemma_A'} by
restricting $H< q^{\frac{1}{2} + \frac{1}{4(r - D)}};$ the requirement
$2H/P\le q$ holds for sufficiently large $q$. Then
\[ T(N,H) \ll_{r,\ep} H^{1-1/r} q^{\frac{r+1 - D}{4r(r - D)}+\ep} .\]	
This completes the proof of Theorem \ref{thm_Vin_conj}.

As already noted, Theorems \ref{thm_Vinogradov_1} through
\ref{thm_Vinogradov_3} hold because Conjecture \ref{conj_VMT} is
trivially true for $d=1,2$ and is now known to be true for $d=3$ by
recent results of Wooley \cite{Woo14aX}. For $d\ge 4$ Wooley
\cite{Woo14aX}, \cite{Woo14bX} has proved the following results
towards Conjecture \ref{conj_VMT}:
\begin{prop}\label{prop_Vin}
For $d \geq 4$ and $r \geq d(d-1)$, 
\beq\label{J4}
J_{r,d}(X) \ll_{r,\ep} X^\ep (X^r + X^{2r-D}).
\eeq
For $d \geq 4$ and $D < r< d(d-1)$ then
\beq\label{J3}
 J_{r,d}(X) \ll_r X^{2r -D + \Del},
\eeq
where the order of magnitude of $\Del=\Del(r,d)$ is $O(d)$,
as specified in \cite{Woo14bX}.
\end{prop}
The result (\ref{J4}) immediately implies Theorem \ref{thm_Vinogradov_3_large}. 
Theorem \ref{thm_Vinogradov_med} follows from applying (\ref{J3})  in
Proposition \ref{prop_S4} to deduce that
\[ \sum_\al \sum_{m=1}^q T(\al;m,\tau)^{2r} \ll_{r,\ep}
Q^D \left(\tau^{2r-D + \Del}q^{1/2} + \tau^r q\right)q^\ep . \]
The argument then proceeds as before, after choosing
$P$ such that 
\[ \frac{1}{2}H q^{- \frac{1}{2(r -D+ \Del)}} \leq P < H q^{- \frac{1}{2(r - D+ \Del)}}\]
in place of (\ref{P_choice}). 

\subsection{A note on $\del$}\label{sec_delta_Vin}
We remark that if $H=q^{1/4+\kappa}$ for some small $\kappa>0$ then
Theorem \ref{thm_Vin_conj} would give a nontrivial bound $Hq^{-\del}$ where 
\[ \del =  \frac{4 \kappa (r-D) - 1}{4r(r-D)} .\]
As a function of $r$ this attains a maximum at the real value
\[ r_{\kappa,d} := D+ \frac{1+\sqrt{4D\kappa +1}}{4 \kappa}.\]
We choose $r$ to be an integer $r=r_{\kappa,d} + \theta$ with $-1/2< \theta \leq 1/2$, and for this choice, $\del$ is approximately
\[ \del = \left( \frac{2\kappa}{1+\sqrt{1+4D\kappa}} \right)^2.\]
For any fixed $d$, as $\kappa \maps 0$, this behaves like 
\[ \del = \kappa^2,\]
which we note is independent of $d$.

\section{Proof of Lemma \ref{lemma_A'}}\label{sec_lemma_A'}
 This is merely a generalization of the proof in Section 4 of
 \cite{HB12}.  The first property in Lemma \ref{lemma_A'} is a direct
 result of the definition of $\Acal(m)$, on using our assumption that
 $0\le N\le q$. 

For the second property we first note that $\Acal(m)\le\Acal(m)^2$
since $\Acal(m)$ is a non-negative integer.  It follows that $S_1\le S_2$.

We now observe that $\Acal(m)^2 $ counts
quadruples $(p,p',a,a')$ for which
\[m \leq \frac{N-aq}{p}  < m +H/P, \;\;\; m \leq \frac{N-a'q}{p'}  < m +H/P.\]
For such a quadruple we must have
\[\left| \frac{N-aq}{p}  - \frac{N-a'q}{p'} \right| \leq H/P.\]
Under this condition there are $O(H/P)$ corresponding values of $m$.
It follows that
\begin{eqnarray}
 \sum_{m} \Acal(m)^2 
	& \ll & HP^{-1} \# \{p, p',a,a' : 0 \leq \left| \frac{N-aq}{p}  - \frac{N-a'q}{p'} \right| \leq H/P \} \nonumber \\
	& \ll &  H P^{-1}  \sum_{p, p' \in \Pscr} \Mcal(p,p'), \label{HPM}
	\end{eqnarray}
 where 
 \[ \Mcal(p,p') =  \# \{a\modd{p}, a' \modd{p'} : 0 \leq \left| \frac{N-aq}{p}  - \frac{N-a'q}{p'} \right| \leq H/P \}.\]
 First consider the case $p=p'$. Then
 \[|a-a'|\le\frac{Hp}{Pq}\le\frac{2H}{q}<1,\]
since $H=o(q)$ in all our theorems.  Thus $a=a'$
so that $\Mcal(p,p) \ll P$ and hence $\sum_{p=p' \in \Pscr}\Mcal(p,p')
\ll P^2$, which makes an satisfactory contribution to (\ref{HPM}).

Next, consider the case $p \neq p'$. 
We choose (by Bertrand's postulate) a prime $l$ such that 
\[ \frac{q}{H} < l \leq \frac{2q}{H}.\]
(Here we use the fact that $H< q$ for large enough $q$.)
Let $M = \left[ \frac{N l}{q} \right]$ or $1+\left[ \frac{N l}{q}
\right]$ be chosen so that $l\nmid M$. Then $|Nl/q - M | \leq 1$ implies that $|N- qM/l| \leq q/l$, so that 
\[ \left| \frac{qM/l - aq}{p} - \frac{qM/l  - a'q}{p'} \right| \leq
\frac{H}{P} + \frac{q}{lp} + \frac{q}{lp'}  \]
for every pair $a,a'$ counted by $\Mcal(p,p')$.  Thus
\[\left|M(p'-p)-(ap'-a'p)l\right| \leq
\frac{pp'Hl}{qP}+p'+p\leq \frac{2pp'}{P}+p'+p\leq 12P.\]
For a given $\delta$ there is at most one way to choose $a,a'$ with
$0\leq a<p$ and $0\leq a'<p'$ which satisfy  $ap'-a'p=\del$. 
 Thus 
\[
  \sum_{p \neq p' \in \Pscr} \Mcal(p,p') \ll \#\{ p \neq p' \in \Pscr, 
|m| \leq 12P: \\   M(p'-p) \con m \modd{l} \}.\]
We chose $M$ so that $l\nmid M$, and hence the condition $M(p'-p) \con
m \modd{l}$ determines $p'-p$ uniquely modulo $l$. Since  by hypothesis $P<q/H<l$
this suffices to determine at most two values for $p'-p$ in $\Z$. So
we may choose $p$ freely and there are then at most two possibilities for $p'$. As a result, after counting up the possible choices for $m$, we conclude that 
   \[ \sum_{p \neq p' \in \Pscr} \Mcal(p,p') \ll P^2.\]
Applying this in (\ref{HPM}), we conclude that
\[  \sum_{m} \Acal(m)^2 \ll H P,\]
as required.

\section*{Acknowledgements}

Pierce was partially supported during this work by a Marie Curie
Fellowship funded by the European Commission and the National Science Foundation on grant DMS-0902658.

\bibliographystyle{amsplain}
\bibliography{NoThBibliography}

\end{document}